\documentclass{amsart}

\usepackage{amsmath,amsthm,amsfonts,amscd, amssymb} 
\usepackage{pictexwd,dcpic}
\usepackage[dvips]{graphicx}

\newtheorem{theorem}{Theorem}[section]
\newtheorem{proposition}[theorem]{Proposition}
\newtheorem{lemma}[theorem]{Lemma}
\newtheorem{corollary}[theorem]{Corollary}
\newtheorem*{main}{Theorem 0.1}
\newtheorem{defthe}[theorem]{Definition/Theorem}

\theoremstyle{definition}
\newtheorem{definition}[theorem]{Definition}

\numberwithin{equation}{section}

\newcommand{\bbC}{{\mathbb C}}
\newcommand{\bbZ}{{\mathbb Z}}

\newcommand{\bbQ}{{\mathbb Q}}
\newcommand{\bbP}{{\mathbb P}}

\newcommand{\bbG}{{\mathbb G}}
\newcommand{\calL}{{\mathcal L}}
\newcommand{\calM}{{\mathcal M}}
\newcommand{\calP}{{\mathcal P}}
\newcommand{\calC}{{\mathcal C}}
\newcommand{\calO}{{\mathcal O}}
\newcommand{\AM}{\overline{M}_{0,0}(\bbP^r,d)}
\newcommand{\AG}{\overline{M}_{0,0}(\bbP^r \times \bbP^1, (d,1))}
\newcommand{\AMN}{\overline{M}_{0,n}(\bbP^r,d)}
\newcommand{\AGN}{\overline{M}_{0,n}(\bbP^r \times \bbP^1, (d,1))}

\DeclareMathOperator{\Pic}{Pic}

\begin{document}
\title[A GIT construction of $\AMN$]{An elementary GIT construction of the moduli space of stable maps}
\author{Adam E. Parker}
\address{Department of Mathematics and Computer Science\\
Wittenberg University\\Springfield, Ohio 45504}
\email{aparker@wittenberg.edu}

\begin{abstract}
This paper provides an elementary construction of the moduli space of stable maps $\overline{M}_{0,0}(\mathbb{P}^r,d)$ as a sequence of ``weighted blow-ups along regular embeddings" of a projective variety. This is a corollary to a more general GIT construction of $\overline{M}_{0,n}(\mathbb{P}^r,d)$ that places stable maps, the Fulton-MacPherson space $\mathbb{P}^1[n]$, and curves $\overline{M}_{0,n}$ into a single context.  
\end{abstract}

\maketitle
Given a projective space $\bbP^r$ and a class $d \in A_1(\bbP^r)\cong \bbZ$, a \emph{n-pointed, stable map of degree $d$} consists of the data $\{\mu:C \to \bbP^r, \{p_i\}_{i=1}^n \}$ where:
\begin{itemize}
	\item $C$ is a complex, projective, connected, reduced, n-pointed, genus $0$ curve with at worst nodal singularities.
	\item $\{ p_i \}$ are smooth points of $C$.
	\item $\mu: C \to \bbP^r$ is a morphism.
	\item $\mu_{*}[C] = d l$, where $l$ is a line generator of $A_1(\bbP^r)$.
	\item If $\mu$ collapses a component $E$ of $C$ to a point, then $E$ must contain at least three special points (nodes or marked points).
	\end{itemize}

We say that two stable maps are \emph{isomorphic} if there is an isomorphism of the pointed domain curves $f:C \to C'$ that commutes with the morphisms to $\bbP^r$.  Then there is a projective coarse moduli space $\AMN$ that parametrizes stable maps up to isomorphism \cite{FP}.  The open locus $M_{0,n}(\bbP^r,d)$ corresponds to maps with a smooth domain while  the \emph{boundary} is naturally broken into divisors $D(N_1, N_2, d_1, d_2)$ where $N_1 \cup N_2$ is a partition of $\{1,2, \dots, n\}$ and  $d_1 + d_2 = d$.  This corresponds to maps where the domain curve has two components, one of degree $d_1$ with the points of $N_1$ on it.

Similarly, we can define stable maps to $\bbP^r \times \bbP^1$ of bi-degree $(d,1)$, and look at the corresponding coarse moduli space $\AGN$.  The boundry again is broken into divisors.  When no confusion is possible, we write $D(N_1, N_2, d_1,d_2)$ where we should instead use $D(N_1, N_2, (d_1,1),(d_2,0))$.

In \cite{P1}, Pandharipande constructs the open $M_{0,0}(\bbP^r,d) \subset \AM$ as the GIT quotient of the open basepoint free locus $U(1,r,d) \subset \oplus^r_0 H^0(\bbP^1,\calO(d))$.  We have a similar construction for the open pointed locus $M_{0,n}(\bbP^r \times \bbP^1, (d,1)) \\ \subset \AGN$.  Our main result is to construct the compact $\AMN$ as a geometric quotient of $\AGN$ by $G=Aut(\bbP^1)$.  

\begin{main}\label{t1}
Let $E$ be an effective divisor such that $-E$ is $\phi$-ample.   Take a linearized line bundle $\calL \in \Pic^G((\bbP^1)^n \times \bbP^r_d)$
such that 
\[
((\bbP^1)^n \times \bbP^r_d)^{ss}(\calL) = ((\bbP^1)^n \times \bbP^r_d)^s(\calL) \neq \emptyset.
\]
Then for each sufficiently small $\epsilon > 0$, the line bundle $\calL' = \phi^*(\calL)(-\epsilon E)$ is ample and 
\begin{multline*}
(\AGN)^{ss}(\calL') = (\AGN)^{s}(\calL')  \\
= \phi^{-1}\{((\bbP^1)^n \times \bbP^r_d)^{ss}(\calL)\}.
\end{multline*}
There is a canonical identification
\[
(\AGN)^s(\calL') /G = \AMN
\]
and a commutative diagram
\[
\begin{CD}
(\AGN)^{s}(\calL') @>f>> \AMN \\
@V{\phi}VV @V{\bar{\phi}}VV \\
((\bbP^1)^n \times \bbP^r_d)^{s}(\calL) @>>> ((\bbP^1)^n \times \bbP^r_d)^{s}(\calL) / G \\
\end{CD}
\]
where  $\bbP^r_d := \bbP((H^0(\bbP^1,\calO(d))^{r+1})$ and $\phi: \AGN \to (\bbP^1)^n \times \bbP^r_d$ is the Givental contraction map \cite{G}.
\end{main}

The eventual goal would be to construct  $\AM$ as sequence of blow-ups of some projective variety.  One benefit of such a construction is the ability to compute the Chow ring of $\AM$, as Keel's Theorem 1 from the appendix of \cite{K} gives the Chow ring of a blow up.    This can't happen.  First, $\AM$ is not smooth.  It has singularities at points corresponding to maps with nontrivial automorphisms.  However, $\AM$ is actually smooth when considered as a stack, and so at best we may hope for a stack analogue of a sequence of blow ups mentioned above.   The second issue seems more serious.  There are no known maps from $\AM$ to anything nice, and a birational map from $\AM$ is exactly what is needed to carry out the above project.  

As corollaries to our GIT construction, we are able to construct a birational map $\bar{\phi}$ from $\AM$. Recently, progress has been made on understanding $\AG$.  For example, in \cite{M}, the above $\phi$ is factored into a sequence of intermediate moduli spaces such that the map between two successive spaces is a ``weighted blow up of a regular local embedding".  As a corollary of the above theorem, we take the quotient of these intermediate spaces and factor $\bar{\phi}$.

In Section 1, we collect some preliminary results and definitions that will be used throughout the paper.  Section 2 identifies the stable locus in $(\bbP^1)^n \times \bbP^r_d$ and explains how to pull it back to $\AGN$.  In Section 3, we prove the above theorem.  Finally, in Section 4 we explain the factorization of $\varphi$ from \cite{M} and construct the intermediate spaces the induced quotient map.  All work is done over $\bbC$.

This paper is part of my thesis written at the University of Texas at Austin under the direction of Prof. Sean Keel.   I wish to extend my sincere gratitude for all of his encouragement, guidance, and patience. 

\section{Preliminaries}
Suppose that we wanted to compactify the space of $n$-pointed degree-$d$ morphisms from $\bbP^1 \to \bbP^r$.   Perhaps after the above discussion, one would expect that $\AMN$ correctly compactifies these objects.  However, since we quotient out by isomorphisms, we only get degree $d$ \emph{un}-parametrized pointed morphisms to $\bbP^r$.  Here we discuss two spaces that do correctly answer this question.

\subsection{Linear Sigma Model}

On one hand, an $n$-pointed, degree-$d$ morphism  $f$ is given by $(r+1)$  homogenous degree $d$ polynomials in two variables, along with a choice of $n$ distinct points on the domain $\bbP^1$.    In the notation of \cite{P1}, these maps correspond to the basepoint free locus
\[
\begin{aligned}
((\bbP^1)^n \setminus \Delta) \times \bbP(U(1, r, d)) &\subset (\bbP^1)^n \times \bbP( \bigoplus_0^r H^0(\bbP^1, \calO_{\bbP^1}(d))) \\
& := (\bbP^1)^n \times \bbP^r_d.
\end{aligned}
\]

Once we pick coordinates on $\bbP^1$, we can consider a closed point on the basepoint free locus as  
\[
[x_1:y_1] \times \dots \times [x_n: y_n] \times [ f_0(x,y) : f_1(x,y) : \dots : f_r(x,y)]
\]
where $[x_s:y_s] \neq [x_t:y_t]$, the $f_j$ don't have any common roots, and scaling doesn't change the map.  The coefficients of these $f_j$ determine a point in projective space $\bbP^r_d := \bbP^{(r+1)(d+1)-1}$.  We will sometimes write $a_i^j$ for the coefficient $x^{d-i}y^i$ on $f_j$ (after choosing the obvious coordinates on $\bbP^r_d.)$   We thus have a simple compactification by allowing the $r+1$ forms to have common roots, and allowing the $n$ points to come together.  This space is sometimes referred to as the linear sigma model.

Moreover, there is a $G$ action on this space, similar to the action examined in \cite{MFK} on binary quantics.  On closed points, the action is given by
	\[
	G \times ((\bbP^1)^n \times \bbP^r_d)\to (\bbP^1)^n \times \bbP^r_d
	\]
\begin{multline*}
	g \cdot [ [x_1:y_1] \times \dots \times [x_n:y_n] \times  f_0(x:y) , \dots , f_r( x:y)] = \\
	 [g[x_1: y_1] \times \dots \times g[x_n:y_n] \times f_0 \circ g^{-1}(x:y),  \dots :f_r \circ g^{-1}(x:y)]
\end{multline*}	
where $g$ and $g^{-1}$ act on $[x:y]$ by matrix multiplication.

\subsection{The Graph Space} There is another, less simple (non-linear) compactification of $((\bbP^1)^n \setminus \Delta) \times \bbP(U(1, r, d))$.  It is clear that this set equals $M_{0,n}(\bbP^r \times \bbP^1, (d,1))$, and thus $\AGN$ provides another compacitification. 

We will refer to the domain curve $C$ for a map in $\AGN$ as a \emph{comb}.  There is an obvious distinguished component $C_0$ on which $\mu|_{C_0}$ will be of degree $(d',1)$.  We will call this component the \emph{handle}.  The other components fit into \emph{teeth} $T_i$, which are (perhaps reducible) genus-$0$, $n_i$-pointed curves meeting $C_0$ at unique points $q_i$.  There is always a representative of the map so that degree $1$ part of  $\mu$ restricted to the handle is the identity.  

The action on $\AGN$ is induced by the action on the image $\bbP^1$.  Namely we have 
	\[
	G \times \AGN \to \AGN
	\]
	\[
	g \cdot [\mu_1 \times \mu_2 : C \to \bbP^r \times \bbP^1] \to [\mu_1 \times g \circ \mu_2 : C \to \bbP^r \times \bbP^1]
	\]

\subsection{The Givental Map}
Recall in \cite{G}, that Givental constructs a projective morphism that relates the graph space and the linear sigma model.

\begin{theorem} (Givental) \label{g} There is a projective morphism
\[
\varphi: \AG \to \bbP^r_d
\]
\end{theorem}
Set theoretically,  consider a point in $\AG$.  As mentioned above, there is a representative 
\[
[\mu:C \to \bbP^r \times \bbP^1 \mbox{ of bi-degree } (d,1) ]
\]
and a component $C_0 \subset C$ such that  $\mu|_{C_0}$ is the graph of  $r+1$ degree $d'$ polynomials $(f_0, \dots, f_r)$ with no common zero.  On the teeth $T_1, \dots, T_s$, $\mu$ has degree $(d_i, 0)$ respectively, and $d_1 + \dots + d_s = d - d'$.  Thus $\mu$ sends $T_i$ into $\bbP^r \times z_i \subset \bbP^r \times \bbP^1$.  Let $h$ be a degree $d-d'$ form  that vanishes at each $z_i$ with multiplicity $d_i$.  Then
\[
\varphi(\mu) = [f_0 \cdot h , f_1 \cdot h, \dots, f_r \cdot h] \in \bbP^r_d
\]
where we read off the coefficients to obtain the point in projective space.  The projective morphism that we consider is thus the product of $\varphi$ with the $n$ evaluation morphisms $ev_i: \AMN \to \bbP^r \times \bbP^1 \to \bbP^1$.  We call it $\phi$.

On $M_{0,n}(\bbP^r \times \bbP^1, (d,1))$, $\phi$ gives the isomorphism with $((\bbP^1)^n \setminus \Delta) \times \bbP(U(1, r, d))$ mentioned above.  

The following lemma is needed when we take the quotients.
\begin{lemma}  The above map
\[
\phi: \AGN \to (\bbP^1)^n \times \bbP^r_d
\]
is equivariant with respect to the above $G$ actions.
\end{lemma}
\begin{proof}We show that both the evaluation morphisms $ev_i$ and the Givental $\varphi$ map are equivariant.  Then their product is as well.

Take a point in $\AGN$.  Choose a representative 
$(\mu:C \to \bbP^r \times \bbP^1, \{p_i \})$.  Write $C = C_0 \cup T_i$ as a comb,  such that $T_i \cap C_0 = q_i$.  Also write $\mu_i = \pi_1 \circ \mu|_{T_i}$.

If we look at the image of the above map under $\varphi$, we see that $\varphi(\mu)$ will be the product of  $r+1$ forms $(f_0, \cdots, f_r)$ of degree $d'$ representing the handle and a form $h$ of degree $d-d'$ that vanishes at the $q_i$ with the correct degrees.  We see that $g\cdot \varphi(\mu)$ will be the product of $(f_0 \cdot g^{-1}, \dots, f_r \cdot g^{-1})$ which are $r+1$ forms of degree $d'$  with no common zero with a form $h'$ of degree $d- d'$ that vanishes at $g(q_i)$ with the same degree that $h$ vanished at $q_i$.

We now need to calculate $\varphi(g\cdot \mu)$.  With the above notation, we see that  $g \cdot (\mu:C \to \bbP^r \times \bbP^1)$ would send $p \in T_i$ to $(\mu_i(p),g(q_i))$, and $p \in C_0$ to $(\mu_0(p),g(p))$.  We  find a representative of this new map that has the degree $1$ part be the identity.

Take the curve $C' = g(C_0) \cup T_i$, where now the teeth $T_i$ are glued to $C_0$ at $g(q_i)$.  Define the map from $C' \to \bbP^r \times \bbP^1$ as $(\mu_0 \circ g^{-1}, id)$ on $g(C_0)$, and agrees with $\mu_i$ on the other teeth.  The corresponding map is isomorphic to $g \cdot \mu$.

We look at the image under the Givental map.  The image will be the $r+1$ degree $d'$ forms $\mu_0 \circ g^{-1}$, along with a form $h$ that vanishes at $g(q_i)$ of the correct degree.  This is the same as $g \cdot \varphi(\mu)$.  This shows that $\varphi$ is equivariant.

That the evaluation morphisms are equivariant is immediate.
\end{proof} 

\subsection{The Forgetful Morphism}
The second map that we will be interested in is the ``forgetful" morphism
\[
f: \AGN \to \AMN
\]
defined by forgeting the map to $\bbP^1$ and collapsing any components that become unstable.   Moreover, since the $G$ action on $\AMN$ is trivial, we automatically have $f$ is $G$ equivariant.

\section{Calculations on $\bbP^r_d$}
Immediately, one would expect that $\AMN$ is the quotient of $\AGN$ by $G$ as $G$ ``takes into account" the map to $\bbP^1$.  The question is, ``How to  take the quotient?"  We will use Geometric Invariant Theory in order to find an open set in $\AGN$ such that the quotient is $\AMN$.

All background concerning GIT will be taken from \cite{D} and \cite{MFK}, though we recall the main theorem here for reference.

\begin{theorem}  \cite{D} \cite{MFK} \label{D} Let $X$ be an algebraic variety, $\calL$ a G-linearized line bundle on $X$.  Then there are open sets $X^{s}(\calL) \subset X^{ss}(\calL)\subset X$ (the \emph{stable and semi-stable loci}), such that the quotient
\[
\pi:X^{ss}(\calL) \to X^{ss}(\calL) /G
\]
is quasi-projective and a ``good categorical quotient".  This says (among other things) that for any other $G$-invariant morphism $g:X^{ss}(\calL) \to Z$, there is a unique morphism $h: X^{ss}(\calL) //G \to Z$ satisfying $h\circ \pi = g$.  If we restrict to $X^{s}(\calL)$ then we have a ``geometric quotient".  This says (among other things) that the geometric fibers are orbits of the geometric points of $X$, and the regular functions on $X^{s}(\calL)/G$ are $G$ -equivariant functions on $X$.
\end{theorem}

Of special interest to us will be when $X$ is proper over $\bbC$ (as in this case) and $\calL$ is ample, as $X^{ss}(\calL) /G$ will be projective \cite{MFK}.

We start with considering the $G$ action on $(\bbP^1)^n \times \bbP^r_d$.  In our case since $G $  acts on a normal, irreduciable, proper variety $X$ (such as $\bbP^r)$.  Then any line bundle admits a unique $G$ linearization.

\begin{proposition}
\[
\Pic^{G}((\bbP^1)^n \times \bbP^r_d) \cong \bbZ^{n+1}
\]
\end{proposition}
\begin{proof}
For any vector $\vec{k}= (k_1, \cdots, k_n,k_{n+1}) \in \bbZ^{n+1}$, we define a line bundle on $(\bbP^1)^n \times \bbP^r_d$ by
\[
\calL_{\vec{k}}=\bigotimes_{i=1}^{n+1} \pi_i^{*}(\calO(k_i))
\]
where $\pi_i$ is projection onto the $i$-th component.  Every line bundle on $(\bbP^1)^n \times \bbP^r_d$ is isomorphic to $\calL_{\vec{k}}$ for some choice of $\vec{k}$ (\cite{hart}).  We need only show that each of these line bundles has one (and only one) linearization.  However, since each $\pi_i$ is $G$-equivariant, and each of the restrictions of $\calL_{\vec{k}}$ to a factor has a unique linearization (\cite{D}), $\calL_{\vec{k}}$ has a canonical $G$ - linearization.
\end{proof}

\begin{corollary}
\[
\calL_{\vec{k}} \text{ is ample } \iff k_i >0
\]
\end{corollary}
\begin{proof}
If all $k_i >0$ then $\calL_{\vec{k}}$ defines the projective embedding
\begin{multline*}
\begin{CD}
(\bbP^1)^n \times \bbP^r_d@>{\text{Veronese}}>> \prod_{i=1}^n \bbP^{(1+k_i)
-1} \times \bbP^{\binom{rd+r+d + k_{n+1}}{rd+r+d}-1} 
\end{CD}\\
\begin{CD}
@>{\text{Segre}}>> 
\bbP^{\left( ( \prod_{i=1}^n 1+k_i) \times {\binom{rd+r+d + k_{n+1}}{rd+r+d}} \right)-1}
\end{CD}
\end{multline*}
On the other hand if some multiple of $\calL_{\vec{k}}$ defines a closed embedding, restricting it to any factor will be ample.   But this is $\calO_{\bbP^1}(k_i)$ (or $\calO_{\bbP^r_d}(k_{n+1})$) and these are ample iff $k_i, k_{n+1} >0$.
\end{proof}

In order to find the stable and semi-stable loci in $(\bbP^1)^n \times \bbP^r_d$, we will look at the image under the above Veronese / Segre maps.  The main point is the following.
\begin{proposition}
Let $\Omega$ be the composition of the Veronese and Segre maps above.  Then 
\[
((\bbP^1)^n \times \bbP^r_d)^{ss}(\calL_{\vec{k}}) = \Omega^{-1}\left\{\left(\bbP^{\left( ( \prod_{i=1}^n 1+k_i) \times {\binom{rd+r+d + k_{n+1}}{rd+r+d}} \right)-1}\right)^{ss}(\calO(1))\right\}
\]
and similarly for the stable locus.
\end{proposition}
\begin{proof}
Call the image projective space $\bbP^{BIG}$.   First, we show that there is an action of $G$ on $\bbP^{k}$ such that the Veronese map $\bbP^1 \to \bbP^k$ is $G$ equivariant.  We need a reparesentation of $G$ in $GL(k+1)$.  We can explicitly write it out, where we choose $[x,y]$ as coordinates on $\bbP^1$ and the obvious coordinates $[x^k: x^{k-1}y: \dots: y^k]$ on $\bbP^k$.  Namely, 
\[
\rho: G \to GL(k+1)
\]
\[
\left(
\begin{matrix}
a & b \\
c & d
\end{matrix}
\right) \to [a_{i,j}]_{i,j=0}^k
\]
where $a_{i,j}=\sum_{n=0}^j \binom{k-i}{n}\binom{i}{j-n} b^n a^{k-i-n} d^{j-n} c^{i-j+n}$.  This is the coefficient  of $x^{k-j}y^j$ in $(ax+by)^{k-i} (cx+dy)^i$,  and is a homomorphism.  We can define the representation of $G$ into $GL(\binom{rd+r+d + k_{n+1}}{rd+r+d})$ similarly.  We now have representations 
\[
\rho_i: G \to GL(k_i +1) \,\,\,\,\, \text{and} \,\,\,\,\,\, \rho_{n+1}: G\to GL\left(\binom{rd+r+d + k_{n+1}}{rd+r+d}\right).
\]
We define the action on $\bbP^{BIG}$ by taking the tensor representation.   This extends to an action on all of $\bbP^{BIG}$.  Thus $\Omega$  is $G$-invariant by construction.

Take the composition $(\bbP^1)^n \times \bbP^r_d \to \Omega((\bbP^1)^n \times \bbP^r_d) \hookrightarrow \bbP^{BIG}$.  We apply the following theorem of \cite{MFK} to each of these arrows.  
\begin{theorem} \cite{MFK} (pg 46)
Assume that $f:X \to Y$ is finite, $G$- equivariant with respect to actions of $G$ on $X$ and $Y$.  If $X$ is proper over $k$ ($\bbC$ for us) and $M$ is ample on $Y$, then
\[
X^{ss}(f^*M) = f^{-1}\{Y^{ss}(M)\}
\]
and the same result holds for the stable locus.
\end{theorem}
Finally, that $\Omega^* \calO(1) = \calL_{\vec{k}}$ is obvious.
\end{proof}
We are now able to determine the stable and semi-stable locus in the linear-sigma model.

\begin{theorem}
Let $\vec{k}=(k_1, k_2, \dots, k_{n+1}) \in \bbZ^{n+1}_+$.  Then $[x_1: y_1] \times \dots \times [x_n: y_n] \times [a_i^j] \in ((\bbP^1)^n \times \bbP^r_d)^{ss}(\calL_{\vec{k}})$ (respectively $\in ((\bbP^1)^n \times \bbP^r_d)^{s}(\calL_{\vec{k}})$) if for every point $p\in \bbP^1$
\[
\sum_{i\,|\, [x_i: y_i] = p} k_i +  d_p \cdot k_{n+1} \leq \frac{1}{2} \left(\sum_{i=1}^n k_i + d\cdot k_{n+1} \right)
\]
(respectively strict inequality holds) where $d_p$ is the degree of common vanishing of the forms $f_0, \dots, f_r$ at $p \in \bbP^1$.
\end{theorem}
\begin{proof}
We prove the theorem by first looking at the action of a maximal torus acting on $\bbP^{BIG}$.  Here, there is only one line bundle, so everything is canonical.  Then we pull back to find the corresponding locus in $(\bbP^1)^n \times \bbP^r_d$.  We then move onto the entire group G.

Let T be the maximal torus of $SL_2(\bbC)$, equal to the image of the 1-parameter subgroup 
\begin{equation*}
\lambda(t) =
\left(
\begin{matrix}
t^{-1} & 0 \\
0 & t
\end{matrix}
\right).
\end{equation*}

We choose coordinates $a^j_i$ on $\bbP^r_d$, where  $a^j_i$ is the coefficient of $x^{d-i}y^i$ in $f_j(x,y)$.  Similarly, we choose the following coordinates on $\bbP^{BIG}$.  For $0 \leq s_i \leq k_i$  $(1\leq i \leq n)$, and $v_{ij}$ such that $\sum_{i=0}^d \sum_{j=0}^r v_{ij} = k_{n+1}$, we have the coordinate $x_i^{k_i - s_i} y_i^{s_i}(a_i^j)^{v_{ij}}$.  Then $T$ acts on $\bbP^{BIG}$ by
\[
\lambda(t) \cdot (x_i^{k_i - s_i} y_i^{s_i}(a_i^j)^{v_{ij}}) \to t^{(\sum_{i=1}^n 2s_i-k_i) + (\sum_{ij} (d-2i)v_{ij})}x_i^{k_i - s_i} y_i^{s_i}(a_i^j)^{v_{ij}}
\]

By the above Lemma, we know that it's enough to compute the semi-stable locus of this action on $\bbP^{BIG}$ and pull it back via the various inclusions and embeddings.  Luckily we know how to compute the semi-stable locus of a torus acting on a projective space.  From \cite{D} we know that a point of projective space is stable (resp semi-stable) with respect to $T$ if and only if $0 \in \mbox{interior} (\overline{wt})$ ( resp $0 \in \overline{wt}$).  In our case, the weight set (wt) is the subset of 
\[
\left\{ -\sum_{i=1}^n k_i - d \cdot k_{n+1}, \dots, \sum_{i=1}^n k_i + d \cdot k_{n+1} \right\}
\]
consisting of powers of $t$ such that the coordinate $x_i^{k_i - s_i} y_i^{s_i}(a_i^j)^{v_{ij}} $ is non zero.  

If the point is unstable, then all the powers $(\sum_{i=1}^n 2s_i-k_i) + (\sum_{ij} (d-2i)v_{ij}) < 0 $ (or perhaps all $>0$.)   So $x_i^{k_i - s_i} y_i^{s_i}(a_i^j)^{v_{ij}} =0$  if  
\begin{multline*}
0 \leq \left(\sum_{i=1}^n (2s_i-k_i) + \sum_{ij} (d-2i)v_{ij} \right)  \iff  \\
0 \leq 2 \sum_{i=1}^n (s_i-k_i) - 2 \sum_{ij} i \cdot v_{ij} + \sum_{i=1}^n k_i + d k_{n+1}  \iff \\
 \sum_{ij} i\cdot v_{ij} + \sum_{i=1}^n (k_i -s_i) \leq \frac{1}{2}\left( \sum_{i=1}^n k_i + d k_{n+1} \right) .
\end{multline*}

Define the following sets in $(\bbP^1)^n \times \bbP^r_d$:
\begin{multline*}
US = \{ [x_1: y_1] \times \dots \times [x_n: y_n] \times [a_i^j] \,\,| \\
 x_i^{k_i - s_i} y_i^{s_i}(a_i^j)^{v_{ij}}= 0\,\, \text{ if } \,\,\ \sum_{ij} i\cdot v_{ij} + \sum_{i=1}^n (k_i -s_i) \leq \frac{1}{2}\left( \sum_{i=1}^n k_i + d k_{n+1}  \right) \}
\end{multline*}
and
\begin{multline*}
X = \{[x_1: y_1] \times \dots \times [x_n: y_n] \times [a_i^j] \,\, | \\
 \frac{1}{2} \left( \sum_{i=1}^n k_i + d k_{n+1} \right)  < \sum_{[x_i:y_i] = [1:0]} k_i + k_{n+1} \cdot d_{[1:0] } \}
\end{multline*}
We show that $US = X$.

First, assume that $X \subset US$.   Let $x  = \{ [x_1: y_1] \times \dots \times [x_n: y_n] \times [a_i^j] \}$ be in $US \setminus X$.  So, $ x_i^{k_i - s_i} y_i^{s_i}(a_i^j)^{v_{ij}}= 0$ if  
\[
 \sum_{ij} i\cdot v_{ij} + \sum_{i=1}^n (k_i -s_i) \leq \frac{1}{2}\left( \sum_{i=1}^n k_i + d k_{n+1}  \right).
\]
But we also have
\[
\sum_{[x_i:y_i] = [1:0]} k_i + k_{n+1} \cdot d_{[1:0] } \leq \frac{1}{2} \left( \sum_{i=1}^n k_i + d k_{n+1} \right) .
\]
Then, take $s_i = 0$ if $[x_i: y_i] = [1:0]$.  And at least one of the $a_{d_{[1:0]}}^j \neq 0$.  For that value of $j$, let $v_{ij}= k_{n+1}$.  Then we have
\[
 \sum_{i=0}^n(k_i - s_i) + \sum_{ij} i\cdot v_{ij}  = \sum_{[x_i:y_i] = [1:0]} k_i + k_{n+1} \cdot d_{[1:0] } \leq \frac{1}{2} \left( \sum_{i=1}^n k_i + d k_{n+1} \right).
 \]
 The coordinate $x_i^{k_i - s_i} y_i^{s_i}(a_i^j)^{v_{ij}}\neq 0$ by construction, which says $x \notin US$, a contradiction.  Thus $US \subseteq X$.

Now, assume that $US \subset X$  Take $y =  \{ [x_1: y_1] \times \dots \times [x_n: y_n] \times [a_i^j] \} \in X \setminus US$.  So $ x_i^{k_i - s_i} y_i^{s_i}(a_i^j)^{v_{ij}}\neq 0$, but 
\[
\sum_{ij} i\cdot v_{ij} + \sum_{i=1}^n (k_i -s_i)\} \leq \frac{1}{2}\left( \sum_{i=1}^n k_i + d k_{n+1}  \right) 
\]
Combining with the fact that $y \in X$, we see that 
\[
\sum_{i=1}^n (k_i - s_i) + \sum_{ij} i \cdot v_{ij} < \sum_{[x_i:y_i] = [1:0]} k_i + k_{n+1} \cdot d_{[1:0] }
\]
Then, for all $i$ with $[x_i:y_i] = [1:0]$, we must have $s_i = 0$.  Thus,
\[
\sum_{[x_i:y_i] = [1:0]} k_i  \leq \sum_{i=1}^n(k_i - s_i).
\]
Similarly, since $(a_i^j)^{v_{ij}} \neq 0$, we know that if $(a_i^j) = 0$, then $v_{ij}=0$.  Thus,
\[
 \sum_{j=0}^r \sum_{i=0}^di \cdot v_{ij} = \sum_{j=0}^r \sum_{i = d_{[1:0]}}^d i \cdot v_{ij} \geq d_{[1:0]} \sum_{ij} v_{ij} = d_{[1:0]} \cdot k_{n+1}
 \]
 Combining these gives our contradiction, showing that $X\subseteq US$ as desired.

If we repeat this calculation, except replace the condition $(\sum_{i=1}^n 2s_i-k_i) + (\sum_{ij} (d-2i)v_{ij}) < 0 $ with $(\sum_{i=1}^n 2s_i-k_i) + (\sum_{ij} (d-2i)v_{ij})  > 0$, we get the following Lemma. 

\begin{lemma}\label{l1}
$[x_1: y_1] \times \dots \times [x_n, y_n] \times [a_i^j]$ is unstable with respect to $T$ if 
\[
\sum_{i \, | \,[x_i:y_i] = [1:0]} k_i + k_{n+1}\cdot d_{[1:0]} > \frac{1}{2} \left( \sum_{i=1}^n k_i + d k_{n+1}\right)
\]
or
\[
\sum_{i \, | \,[x_i:y_i] = [0:1]} k_i + k_{n+1} \cdot d_{[0:1]}> \frac{1}{2} \left( \sum_{i=1}^n k_i + d k_{n+1}\right)
\]
\end{lemma}

We are now ready to move onto stability with respect to G.  Suppose that $[x_1: y_1] \times \dots \times [x_n, y_n] \times [a_i^j]$ is stable with respect to G and there is a point $p$ in $\bbP^1$ such that 
\[
\sum_{[x_i:y_i] = p} k_i + k_{n+1} \cdot d_p > \frac{1}{2} \left( \sum_{i=1}^n k_i + d k_{n+1}\right).
\]
Let $g\in G$ map $ p \to [1:0]$.  Then $g\cdot [x_1: y_1] \times \dots \times [x_n, y_n] \times [a_i^j]$ is unstable with respect to T, and $[x_1: y_1] \times \dots \times [x_n, y_n] \times [a_i^j]$ is unstable with respect to G, contradicting the assumption.

Now assume that $[x_1: y_1] \times \dots \times [x_n, y_n] \times [a_i^j]$ is unstable, but has no point $p$ such that 
\[
\sum_{i | [x_i:y_i] = p} k_i + k_{n+1} \cdot d_p > \frac{1}{2} \left( \sum_{i=1}^n k_i + d k_{n+1}\right).
\]
Then, there is some maximal torus $T'$ with which $[x_1: y_1] \times \dots \times [x_n, y_n] \times [a_s]$ is unstable.  For any maximal torus in G, there is $g\in G$ such that $gT'g^{-1} = T$.  Then we have that $g \cdot  [x_1: y_1] \times \dots \times [x_n, y_n] \times [a_i^j]$ is unstable with respect to $T$,  hence  must have either $[1:0]$ or $[0:1]$ satisfiying Lemma \ref{l1}.  Then $[x_1: y_1] \times \dots \times [x_n, y_n] \times [a_i^j]$ has $g^{-1}[1:0]$ satisfying Lemma \ref{l1}.
\end{proof}

We are now ready to describe the chamber decomposition of the ample cone of  $\Pic^G ((\bbP^1)^n \times \bbP^r_d)$.  As a first step we normalize our line bundle so that we form the simplex
\[
 \Delta = \left\{ (k_1, k_2, \dots, k_{n+1}) | \sum_{i=1}^n k_i + d \cdot k_{n+1} = 2 \right\}
\]
Then for each subset $I \in (1,2, \dots, n)$  and each integer $0 \leq d_{I} \leq d$, we get a wall $W_{I, d_I}$ given by 
\[
\sum_{i \in I} k_i + d_I \cdot k_{n+1} =1
\]
and the walls break $\Delta$ into chambers.  Following \cite{H2}, we mention the following obvious statements.

\begin{enumerate}
\item 
\[
W_{S, d_S} = W_{S^c, d-d_S}
\]
\item Each interrior wall divides $\Delta$ into two parts
\[
\left\{
(k_1, k_2, \dots, k_{n+1}) | \sum_{i \in I} k_i + d_I \cdot k_{n+1} \leq 1\right\}
\]
and
\[
\left\{(k_1, k_2, \dots, k_{n+1}) | \sum_{i \in I} k_i + d_I \cdot k_{n+1} \geq 1\right\}
\]
\item  Two vectors $\vec{k}= (k_1, \dots, k_{n+1})$ and $\vec{k'} =  (k'_1, \dots, k'_{n+1})$ lie in the same chamber if for all $I \subset \{ 1, 2, \dots \}$ and $0 \leq d_I \leq d$ then
\[
\sum_{i \in I} k_i + d_I \cdot k_{n+1} \leq 1 \iff  \sum_{i \in I}k' _i + d_I \cdot k'_{n+1} \leq 1
\]
This means that vectors in the same chamber will define the same stable and semi-stable loci, and hence the same quotient.
\item There are semi-stable points that aren't stable iff  $\vec{k}$ lies on a wall.
\end{enumerate}

Recall that our goal isn't to take the GIT quotient of $(\bbP^1)^n \times \bbP^r_d$, but of $\AGN$, and at this point we haven't said anything about the stable or semi-stable loci in $\AGN$.  We are able to pull back the stable locus via $\phi$, by the following Theorem of Yi Hu.
\begin{theorem} \cite{H1}\label{H1}
Let $\pi:Y \to X$ be a $G$-equivariant projective morphism between two (possibly singular) quasi-projective varieties.  Given any linearized ample line bundle $L$ on $X$, choose a relatively ample linearized line bundle $M$ on $Y$.  Assume moreover that $X^{ss}(L)=X^{s}(L)$.  Then there exists a $n_0$ such that when $n\geq n_0$, we have 
\[
Y^{ss}(\pi^* L^n \otimes M) = Y^{s}(\pi^* L^n \otimes M) = \pi^{-1}\{X^s(L)\}
\]
\end{theorem}

\begin{figure}[h!]\label{f7}
\centering
\includegraphics[scale=.4]{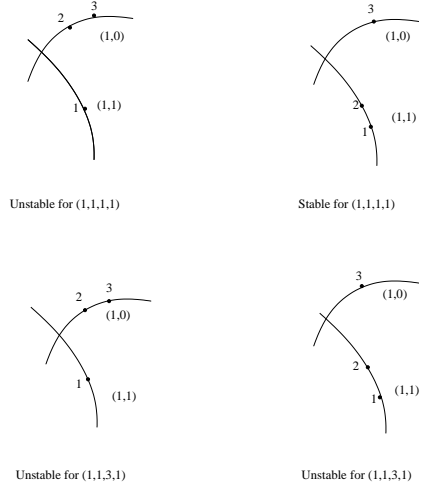}
\caption{Stable Locus of $\overline{M}_{0,3}(\bbP^r \times \bbP^1, (2,1))$}
\end{figure}

For example, the locus of maps in $\AG$ that are stable will be maps such that no tooth of the comb $C$ has degree $\geq d/2$.  See Figure \ref{f7} to see how the stable locus depends on the line bundle.

\section{The Geometric Quotient}

We are now ready to  present our GIT description of $\AMN$.  The construction is similar to that of $\overline{M}_{0,n}$ from \cite{HK}.  We state it similarly.  First let $E$ be an effective divisor with support the full exceptional locus of $\phi$, such that $-E$ is $\phi$ ample.  Such an $E$ exists by the following Lemma from \cite{mori}.
\begin{lemma} \cite{mori} (pg. 70) \label{mori}
Let $f:X \to Y$ be a birational morphism.  Assume that $Y$ is projective and $X$ is $\bbQ$-factorial.  Then there is an effective $f$-exceptional divisor $E$ such that $-E$ is $f$-ample.
\end{lemma}
 $\AGN$ is $\bbQ$-factorial because it is locally the quotient of a smooth scheme by a finite group.  
 
\begin{main} For each linearized line bundle $\calL \in \Pic^G((\bbP^1)^n \times \bbP^r_d)$
such that 
\[
((\bbP^1)^n \times \bbP^r_d)^{ss}(\calL) = ((\bbP^1)^n \times \bbP^r_d)^s(\calL) \neq \emptyset
\]
and for each sufficiently small $\epsilon > 0$, the line bundle $\calL' = \phi^*(\calL)(-\epsilon E)$ is ample and 
\begin{multline*}
(\AGN)^{ss}(\calL') = (\AGN)^{s}(\calL') \\
= \phi^{-1}\{((\bbP^1)^n \times \bbP^r_d)^{ss}(\calL)\}.
\end{multline*}
There is a canonical identification
\[
(\AGN)^s(\calL') /G = \AMN
\]
and a commutative diagram
\[
\begin{CD}
(\AGN)^{ss}(\calL') @>f>> \AMN \\
@V{\phi}VV @VVV \\
((\bbP^1)^n \times \bbP^r_d)^{ss}(\calL) @>>> ((\bbP^1)^n \times \bbP^r_d)^{ss}(\calL) / G \\
\end{CD}
\]
where $\phi$ is the generalized Givental map, $f$ is the forgetful morphism.
\end{main}
\begin{proof}
For the first two statements, we apply the above Theorem \ref{H1} of Hu.  Following the notation from \cite{HK}, let $U$ be the semi-stable locus in $(\bbP^1)^n \times \bbP^r_d$ for the above action of $G$ corresponding to $\calL_{\vec{k}}$.  Recall that this corresponds to $([x_i, y_i], f_0, \dots f_r)$ such that for any $p \in \bbP^1$, we have 
\[
\sum_{i | [x_i: y_i] = p} k_i + k_{n+1} \cdot d_p \leq \frac{1}{2} \left( \sum_{i=1}^n k_i + k_{n+1} \cdot d \right).
\]
Let $U' = \phi^{-1}(U)$.  Let the corresponding quotients be $Q$ and $Q'$.  We have the obvious composition of $G$  invariant maps:
\[
U' \to U \to Q.
\]
And by the universal properties of GIT quotients, we get a proper birational map $Q' \to Q$.  Similarly, since $G$ acts trivially on $\AMN$, we have by the universal property again a proper birational map from $Q' \to \AMN$.  We will show that this is an isomorphism by showing that both sides have the same Picard number.  This is enough since both sides are $\bbQ$-factorial.  
\[
\rho(Q') = \rho(U') = \rho(U) +e(U) = \rho (Q) + e(U)
\]
where $e(u)$ is the number of $\phi$ exceptional divisors that meet $U'$.  Since $\phi$ is an isomorphism on the open locus $M_{0,n}(\bbP^r \times \bbP^1, (d,1)) \subset \AGN$, we need only look at the boundary divisors in $\AGN$.  We use Lemma \ref{codim} to see which divisors are exceptional.

\begin{lemma}\label{codim}
\[
\phi(D(N_1,N_2,d_1, d_2)) \subset (\bbP^1)^n \times \bbP^r_d
\]
has codimension $|N_2| + (r+1)d_2-1$.
\end{lemma}
\begin{proof}
The idea for this proof comes from Kirwan (\cite{k1}).  First notice that 
\[
\phi(D(N_1, N_2, d_1,d_2)) = (p_1, p_2, \dots, p_n, f_0, \dots, f_r)
\]
where $ p_i =p_j \text { for } i,j \in N_2$, and each of the $f_s$ vanish at that point of multiplicity $d_2$.

First, we calculate the codimension of $(p_1, p_2, \dots, p_n, f_0, \dots, f_r)$ where each of the $p_i$ is $[0:1]$, and where each of $f_j$ has a zero of order $d_2$ at $[0:1]$ and an order of $d-d_2 = d_1$ at $[1:0]$ (i.e. each $f_j$ consists of only the monomial $x^{d_2}y^{d_1}$).  It  is clear this has codimension $n+rd+r+d - r = n+ (r+1)d$.  If we remove the condition that each $f_j$ have a root of order $d_1$ at $[1:0]$, then we allow each $f_j$ to have higher powers of $x$.  We also remove the condition that those $p_i$ with $i\notin N_2$ are equal to $[0:1]$.  Thus we see that the the set of $(p_1, p_2, \dots, p_n, f_0, \dots, f_r)$ such that $[0:1] = p_i$ for $i \in N_2$, and each $f_i$ vanishes at $[0:1]$ with multiplicity $d_2$ has codimension
\[
n +(r+1)d - (r+1) d_1 - |N_1| = |N_2| + (r+1)d_2.
\]
Finally, we act on this set by $G$.  We subtract one from the above codimension because $G$ has dimension two, but we don't count the two dimensional stabilizer of $[0:1]$.
\end{proof}

Next we show that $\rho(Q')$ is independent of the chamber from where $\calL_{\vec{k}}$ comes from.  We check that as we cross a wall $W_{I,d_I}$, $\rho(Q')$ doesn't change.  Let our two open sets be $U_1$ and $U_2$.  Recall that $W_{I,d_I}$ breaks our chamber into two parts
\[
\left\{(k_1, \dots, k_n, k_{n+1}) | \sum_{i \in I} k_i + d_I \cdot k_{n+1} \leq 1 \right\}
\]
and
\[
\left\{(k_1, \dots, k_n, k_{n+1}) | \sum_{i \in I} k_i + d_I \cdot k_{n+1} \geq 1 \right\}
\]
so suppose that $U_1$ meets the first set.  Notice that $U'_1$ and $U'_2$ meet the same divisors $D(N_1, N_2, d_1, d_2)$ except that $U'_1$ meets $D(I^c, I, d-d_I, d_I)$ but not $D(I, I^c, d_I, d-d_I)$.  Similarly, $U'_2$ meets $D(I, I^c, d_I, d-d_I)$ but not $D(I^c, I, d-d_I, d_I)$.  

If $2<|I| , 1< r$ and $1 \leq d_I \leq d$, or $r=1$ and $1<d_i \leq d$, then $Q_1 \dashrightarrow Q_2$ is a small modification (an isomorphism in codimension $1$ in the notation of \cite{mori}).  Hence $\rho(Q_1) = \rho (Q_2)$, and it's clear that $e(U_1) = e(U_2)$.

If $2=|I| \text{ and } d_I = 0$, then we see that $Q_1 \dashrightarrow Q_2$ contracts the divisor $(p_1, \dots, p_n, f_0, \dots f_r)$ where $p_i = p_j, i,j \in I$.  Therefore $\rho(Q_1) = \rho (Q_2)+1$.  However, by Lemma \ref{codim}, we see that the divisor $D(I^c, I, d, 0)$ with $|I|=2$ lying over $U_1$ is not exceptional, while it's complement $D(I, I^c, 0 ,d)$ lying over $U_2$ is exceptional.   Hence $e(U_2) = e(U_1)+1$.  Putting these together we see that $\rho(Q'_1) = \rho(Q'_2)$ as desired. 

If $r=1, |I|=0, d_I =1$, then we see that $Q_1 \dashrightarrow Q_2$  contracts the divisor $(p_1, \dots, p_n, f_0, f_1)$ where $f_0, f_1$ have a common root.  Therefore $\rho(Q_1) = \rho (Q_2)+1$.  However, by Lemma \ref{codim}, we see that the divisor $D(N,0,d-1,1)$ lying over $U_1$ is not exceptional, while it's complement $D(0,N, 1, d-1)$ is contracted.  Hence $e(U_2) = e(U_1)+1$.  Putting these together we see that $\rho(Q'_1) = \rho(Q'_2)$ as desired. 

Finally, we prove the Theorem for one vector of one chamber.  Here we look at all divisors $D(N_1, N_2, d_1,d_2) \subset \AGN$.   We have $2^n$ ways to distribute the $n$ points on the domain curve, and we can label the collapsed component with any degree $\leq d$.  Hence there are $2^n ( d+1) $ potential configurations, however the configurations $D(I, I^c, d, 0)$ are not stable maps if $|I| = n$ or $n-1$.  Hence there are $2^n(d+1) - n-1$ total boundary divisors in $\AGN$.  We need to determine how many are stable (with respect to the group).  We do several calculations, as a given linearization $\vec{k}$ may lie in a maximal chamber for certain values of $d,n$, but lie on a wall for others.  All the calculations are very similar though.  Assume that $r>1$.

\begin{itemize}
\item  CASE 1 ($d+n$ odd.  $d>n$)  We choose the linearization corresponding to $(1,1,1, \dots, 1, 1)$.  We count the unstable divisors, i.e. the number of $D(N_1, N_2, d_1, d_2)$  such that 
\[
|N_2| + d_2 \geq \frac{d+n+1}{2}
\]
Any divisor $D(N_1, N_2, d_1, d_2)$ with $\frac{d+n+1}{2} \leq d_2 \leq d$ is unstable.  There are $2^n(\frac{d-n+1}{2})$ of these.  Thus the total number of unstable divisors is

\[
2^n\left(\frac{d-n+1}{2}\right) + \overbrace{(2^n - \binom{n}{0})}^{\# \text{ with } d_2 = \frac{d+n-1}{2}}+\overbrace{(2^n - \binom{n}{0} - \binom{n}{1})}^{ \# \text{ with } d_2 = \frac{d+n-3}{2}} +
\]
\[
\dots +\overbrace{2^n - \binom{n}{0} - \binom{n}{1} - \dots - \binom{n}{n-1}}^{\# \text{ with } d_2 = \frac{d-n+1}{2}}
\]

\begin{align*}
&= 2^n\left(\frac{d-n+1}{2}\right) +n2^n -n\binom{n}{0} - (n-1)\binom{n}{1} - \dots - 1\binom{n} {n-1} \\
& = 2^n\left(\frac{d-n+1}{2}\right) +n2^n - n \binom{n}{n} - (n-1)\binom{n}{n-1} - \dots - 1\binom{n}{1}\\
& = 2^n\left(\frac{d-n+1}{2}\right) + n2^n - \sum_{i=1}^{n}(i) \binom{n}{i} \\
&= 2^n\left(\frac{d-n+1}{2}\right) +n2^n - n2^{n-1} = 2^{n-1}(d+1)
\end{align*}
Hence the total number of stable divisors is $2^{n-1}(d+1)-1-n$ (stable with respect to the group).  The number which are $\phi$ exceptional are all except those where $I^c = 2$ (by Corollary \ref{codim}).  So there are $2^{n-1}(d+1)-1-n-\binom{n}{2}$ $\phi$-exceptional divisors.  Thus, since $\rho(Q) = n+1$, 
\begin{align*}
\rho(Q') = \rho(Q) + e(U)  &= n+1 + 2^{n-1}(d+1) - 1 -n -\binom{n}{2} \\ 
& = 2^{n-1}(d+1)-  \binom{n}{2}.
\end{align*}

\item  CASE 2 ($d+n$ odd,  $d <n$)  We again choose the linearization corresponding to $(1,1,1, \dots, 1,1)$.  
\item  CASE 3 ($d+n$ even, $n$ odd)  We choose the linearization corresponding to $(1,1, \dots, 1, 2)$.  
\item  CASE 4 ($d+n$ even, $n$ even)  We choose the linearization corresponding to $(1,2,2, \dots, 2, 1)$.  
\end{itemize}

Note that care must be taken when  $n=2$ and $d=1,2$.  For here, $\rho((\bbP^1)^n \times \bbP^r_d) = 2$  for the given linearizations (instead of the expected $3$).  This is because the unstable locus contains a divisor.   However, we wouldn't need to subtract out the divisor $D(0,2,d_1, d_2)$ for not being $\phi$ exceptional, because it would have been unstable with respect to the group.  Thus, the sums work out to be the same.

We have shown for every line bundle such that the stable locus equals the semi-stable locus that $\rho(Q') = 2^{n-1}(d+1) -\binom{n}{2}.$  From \cite{P1} we know
\[
\rho(\AMN) = 2^{n-1}(d+1) - \binom{n}{2}
\]
which completes the proof for $r>1$.

When $r=1$ we repeat the above construction.  Here we see, by Corollary \ref{codim}, that the divisor $D(N, 0, d-1,1)$ is not $\phi$-exceptional.  So we subtract one from the above count of $\phi$-exceptional divisors.  Thus
\[
\rho(Q') = \rho(Q) + e(U) = n+1 + 2^{n-1}(d+1) - 2 -n -\binom{n}{2} = 2^{n-1}(d+1)-  \binom{n}{2} - 1
\]
An immediate consequence of Theorem 4.4 in \cite{BF} gives $\rho(\overline{M}_{0,n}(\bbP^1,d))$ and it agrees with the above calculation.

In the case when $r=0$, then $d=0$ or else the moduli space is empty.  Thus $\bbP^0_0 = pt$.  The calculation follows as now we are only dealing with stable curves, and not stable maps and was proven originally in \cite{HK}.  We have that 
\[
\rho(Q') = 2^{n-1} - \binom{n}{2} -1 = \rho(\overline{M}_{0,n})
\]
\end{proof}

There are three immediate corollaries that are interesting.  The first is a new proof of a result of Keel and Hu in \cite{HK}.  By letting $d,r=0$ in the above Theorem we have.
\begin{corollary} \label{c2} \cite{HK} For each linearized line bundle $\calL \in \Pic^G((\bbP^1)^n)$
such that 
\[
((\bbP^1)^n)^{ss}(\calL) = ((\bbP^1)^n)^s(\calL) \neq \emptyset
\]
and for each sufficiently small $\epsilon > 0$, the line bundle $\calL' = \phi^*(\calL)(-E)$ is ample and 
\[
(\bbP^1[n])^{ss}(\calL') = (\bbP^1[n])^{s}(\calL') = \phi^{-1}((\bbP^1)^n)^ss(\calL)
\]
There is a canonical identification
\[
(\bbP^1[n])^s(\calL') /G = \overline{M}_{0,n}
\]
and a commutative diagram
\[
\begin{CD}
(\bbP^1[n])^{ss}(\calL') @>f>> \overline{M}_{0,n} \\
@V{\phi}VV @VVV \\
((\bbP^1)^n)^{ss}(\calL) @>>> ((\bbP^1)^n )^{ss}(\calL) / G) \\
\end{CD}
\]
\end{corollary}
\begin{proof}  
In the case when $d=1$ and $r=1$, we have that 
\[
\overline{M}_{0,n}(\bbP^1,1) = \bbP^1[n]
\]
where $\bbP^1[n] $ is the Fulton-MacPherson compactification of $n$ points on $\bbP^1$.
The Fulton-MacPherson map $\bbP^1[n] \to (\bbP^1)^n$ is exactly the product of evaluation morphisms $\phi$.
\end{proof}

Secondly, we find that the Grassmannian of lines is a GIT quotient of a projective space.  
\begin{corollary}  \label{c3} The Grassmannian of lines in $\bbP^r$ is the GIT quotient of $\bbP^r_1=\bbP^{2(r+1)-1}$ by the above action of $G$.  
\end{corollary}
\begin{proof}
We know that $\overline{M}_{0,0}(\bbP^r,1) =M_{0,0}(\bbP^r,1) = \bbG(1,r).$ By the Theorem 
\[
(\overline{M}_{0,0}(\bbP^r \times \bbP^1, (1,1)))^{s} / G  =  \bbG(1,r).
\]
But $\overline{M}_{0,0}(\bbP^r \times \bbP^1, (1,1))^{s} = {M}_{0,0}(\bbP^r \times \bbP^1,(1,1)) \cong (\bbP^r_1)^{s}$.  Note that when $n=0$, there is only one ample line bundle (up to multiple) on the linear sigma model and it has a unique linearization. 

\end{proof}

The third corollary constructs $\AM$ as a sequence of intermediate moduli spaces.  It requires more background, and follows.


\section{Intermediate Moduli Spaces}
In the case when $n=0$, we obtain as a corollary a factorization of the induced map
\[
\bar{\varphi}:\AM \to (\bbP^r_d)^{s}/Aut(\bbP^1)
\]
into a sequence of intermediate moduli spaces, such that the map between successive spaces is ``almost" a blow up.  We will use a factorization of the Givental map $\varphi$ presented in \cite{M}.  We give the necessary definitions and results here.

Recall the construction of $\AM$ presented by Fulton and Pandharipande in\cite{FP}.  Given a basis of hyperplanes $\bar{t} \in H^0(\bbP^1, \calO(1))$, there is an open subset $U_{\bar{t}}\subset \AM$ such that if we pull back those hyperplanes, the corresponding domain curves along with the sections will be $(r+1)d$ pointed stable \emph{curves}.   By choosing an ordering on the sections, we get an \'{e}tale rigidification of that open set by a smooth moduli space denoted $M_{0,0}(\bbP^r,d,\bar{t})$ that is a $(\bbC^*)^r$ bundle over $\overline{M}_{0,d(r+1)}$.  
\[
\begin{CD}
\overline{M}_{0,d(r+1)} \supset B@<{(\bbC^*)^{r}}<< {M}_{0,0}(\bbP^r,d, \bar{t}) @>{(S_d)^{r+1}}>> U_{\bar{t}}\subset \AM \\.
\end{CD}
\]
$\AM$ is then constructed by gluing together the $U_{\bar{t}}$ for different choices of $\bar{t}$.

In \cite{M}, Musta\c{t}\v{a} constructs similar rigifications of $\AG$ and $\bbP^r_d$ that are $(\bbC^*)^r$ bundles over $\bbP^1[d(r+1)]$ and $(\bbP^1)^{d(r+1)}$ where $\bbP^1[d(r+1)]$ is the Fulton-MacPherson compactification of $d(r+1)$ points on $\bbP^1$. 
\[
\begin{CD}
\bbP^1[d(r+1)] \supset B@<{(\bbC^*)^{r}}<< G(r,d,\bar{t}) @>{(S_d)^{r+1}}>> U_{\bar{t}}\subset \AG \\
@V{F-M}VV @V{\varphi(\bar{t})}VV @V{\varphi}VV \\
(\bbP^1)^{d(r+1)} @<{(\bbC^*)^{r}}<< \bbP^r_d(\bar{t}) @>{(S_d)^{r+1}}>> \bbP^r_d \\
\end{CD}
\]
The factorization of $\varphi$ is obtained by gluing together the pull back of the following factorization of the Fulton-MacPherson map.
\begin{defthe} (pg 13)
Consider a degree 1 morphism $\phi: C \to \bbP^1$ having as a domain $C$ a rational curve with $N$ marked points.  The morphism will be called $n$ - stable if
\begin{enumerate}
\item Not more than $N-n$ of the marked points coincide.
\item  Any ending curve that is not the parametrized component contains more than $N-n$ points.
\item  All the marked points are smooth, and every component that is not the parametrized component has at least three distinct special points.
\end{enumerate}
There is a smooth projective moduli space $\bbP^1[N,n]$ for families of $n$-stable degree 1 morphisms.  Moreover $\bbP^1[N,n]$ is the blow up of $\bbP^1[N,n-1]$ along the strict transforms of the $n$ - dimensional diagonals in $(\bbP^1)^N$.
 \end{defthe}
There is an analogous factorization of $\varphi(\bar{t})$.
\begin{defthe}(pg 22)
A $(\bar{t},d,k)$ - acceptable family of morphism over $S$ is given by the following data
\[
(\pi:\calC \to S, \phi: \calC \to \bbP^1, \{q_{i,j}\}_{0 \leq i \leq n , 1 \leq j \leq d }, \calL, e)
\]
where
\begin{enumerate}  
\item  The family $(\pi:\calC \to S, \phi: \calC \to \bbP^1, \{q_{i,j}\}_{0 \leq i \leq n , 1 \leq j \leq d })$ is a $(r+1)(k-1) +1$ stable family of degree $1$ morphisms to $\bbP^1$.
\item  $\calL$ is a line bundle on $\calC$
\item  $e:\calO_{\calC}^{r+1} \to \calL$ is a morphism of sheaves with $\pi_*e$ nowhere zero and that, via the natural isomorphism $H^0(\bbP^n,\calO(1)) \cong H^0(S \times \bbP^1, \calO_{S\times \bbP^1}^{r+1})$ we have
\[
(e(\bar{t}_i)=0) = \sum_{j=1}^d q_{i,j}
\]
\end{enumerate}
There is a smooth moduli space $\bbP^r_d(\bar{t},k)$ for these families that is a torus bundle over an open subset of $\bbP^1[(r+1)d, (r+1)(k-1) +1]$.
\end{defthe}

Finally, \cite{M} creates global objects factoring $\varphi$ the same way that $\AM$ was constructed in \cite{FP}.
\begin{defthe}(pg 23)
A $(d,k)$ - acceptable family of morphism is given by the following data
\[
(\pi:\calC \to S, \mu=(\mu_1, \mu_2): C \to \bbP^r \times \bbP^1, \calL, e)
\]
where:
\begin{enumerate}
\item  $\calL$ is a line bundle on $\calC$ which, together with the morphism 
\[
e:\calO_{\calC}^{r+1} \to \calL
\]
determines the rational map $\mu_1:\calC \to \bbP^r$.
\item  For any $s \in S$ and any irreducible component $C'$ of $C_s$, the restriction $e_{C'}: \calO_{\calC'}^{r+1} \to \calL_{C'}$ is non-zero.
\item  For any $s\in S$, $deg \calL_{C_s} = d$ and the image $e_{C_s}(H) \in H^0(\calC_s, \calL_{\calC_s})$ of a generic section $H\in H^0(\calC_s, \calO_{\calC_s}^{r+1})$ determines the structure of a $(r+1)(k-1)+1$ stable morphism on $\mu_2:\calC_s \to \bbP^1$
\end{enumerate}
There is a projective coarse moduli space $\bbP^r_d(k)$ for these objects.
\end{defthe}

It is exactly these objects that we take the quotient of by $Aut(\bbP^1)$.  After we quotient, there is no longer a parametrized component to refer to in the above definitions.  However, when $d$ is odd, there is a unique component of the domain curve that will play this role.

\begin{proposition}
Let $C$ be a connected genus-$0$ curve such that each edge is labeled with a number $d$.  If $\sum d_i$ is odd, then there is a unique irreducible component $\bar{C}$ such that if $C$ is a comb with handle $\bar{C}$, no tooth has sum of degrees $> d/2$
\end{proposition}
\begin{proof}
Let $\{C\}_{d/2}$ be the set of all connected subcurves of degree $\geq d/2$.  Intersect all such subcurves.  There is a unique component in the intersection that will be $\bar{C}$.
\end{proof}

\begin{definition}
A $(d,k)^*$ acceptable morphism is given by the following data.
\[
(\pi:\calC \to S, \mu:\calC \to \bbP^r, \calL, e)
\]
where:
\begin{enumerate}
\item  $\calL$ is a line bundle on $\calC$ which, together with the morphism 
\[
e:\calO_{\calC}^{r+1} \to \calL
\]
determines the rational map $\mu: \calC \to \bbP^r$.
\item  For any $s \in S$ and any irreducible component $C'$ of ${\calC}_s$, the restriction $e_{C'}: \calO_{\calC'}^{r+1} \to \calL_{C'}$ is non-zero.
\item  For any $s\in S$, $deg \calL_{C_s} = d$ and the image $e_{C_s}(H) \in H^0(\calC_s, \calL_{\calC_s})$ of a generic section $H\in H^0(\calC_s, \calO_{\calC_s}^{r+1})$ determines the structure of a $(r+1)(k-1)+1$ - stable rigid morphism where $\bar{C}_{s}$ plays the role of the parametrized component in the definition of a $n$ stable degree $1$ morphism above.
\end{enumerate}
\end{definition}
\begin{corollary}
There is a projective coarse moduli space $\overline{M}_{0,0}(\bbP^r,d,k)$ for families of $(d,k)^*$ acceptable morphisms.
\end{corollary}
\begin{proof}
We show that 
\[
\overline{M}_{0,0}(\bbP^r,d,k) := (\bbP^r_d(k))^{s}/Aut(\bbP^1)
\]
satisfies the properties of a coarse moduli space.  This quotient is constructed identically to that from Theorem 0.1.  We pull back the stable locus from $\bbP^r_d$ by \ref{H1}.  Again, the stable locus in $\bbP^r_d(k)$ will be those $(d,k)$ - acceptable maps such that no tooth has degree $> d/2$.  The universal properties of this space are inherited from the universal properties of $\bbP^r_d(k)$ as well as the universal properties of a categorical quotient.  

First we need to show that there is a natural transformation of functors 
\[
\phi:\overline{\calM}_{0,0}(\bbP^r,d,k) \to Hom_{Sch}(*, \overline{M}_{0,0}(\bbP^r,d,k))
\]
where $\overline{\calM}_{0,0}(\bbP^r,d,k)$ is the obvious moduli functor $\{ schemes \} \to \{ sets \}$.  Given a family of $(d,k)^*$ - acceptable morphisms 
\[
(\pi:\calC \to S, \mu:\calC \to \bbP^r, \calL, e)
\]
we can get a $(d,k)$ - acceptable morphism 
\[
(\pi:\calC \to S, \mu=(\mu_1, \mu_2): C \to \bbP^r \times \bbP^1, \calL, e)
\]
by taking $\mu_2: \calC \to \bbP^1$ to be identity on $\bar{C}_s$ and constant on the other components.  This will lie in the stable locus by construction, and thus gives a map $S \to (\bbP^r_d(k))^{s}$.  Composing with the quotient gives an element of $Hom_{Sch}(S, \overline{M}_{0,0}(\bbP^r,d,k))$.

We need to show that if given a scheme $Z$ and a natrual transformation of functors $\psi: \overline{\calM}_{0,0}(\bbP^r,d,k) \to Hom_{Sch}(*, Z)$,  there exists a unique morphism of schemes
\[
\gamma:\overline{M}_{0,0}(\bbP^r,d,k) \to Z
\]
such that $\psi = \tilde{\gamma} \circ \phi$.   By the above, we have a functor
\[
\overline{\calM}_{0,0}(\bbP^r,d,k) \to (\calP^r_d(k))^s
\]
as such we get a functor
\[
\bar{\psi}: ( \calP^r_d(k))^s \to Hom_{Sch}(*, Z)
\]
which by representability gives a map
\[
\bar{\gamma}: (\bbP^r_d(k))^s \to Z.
\]
This map is $G$ equivariant by construction, hence factors though the quotient
\[
\gamma:\overline{M}_{0,0}(\bbP^r,d,k) \to Z
\]
\end{proof}

We can sum up this corollary with the following figure

\begin{figure}[h!]
\begindc{\commdiag}[24]
\obj(11,1)[B]{$\bbP^r_d /G$}
\obj(1,3)[C]{$(\bbP^1)^N = \bbP^1(N,0) $}
\obj(5,3)[D]{$\bbP^r_d(\bar{t})$}
\obj(1,4)[E]{$\bbP^1[N, 1]\,\,\,\,\,\,\,$}
\obj(5,4)[F]{$\bbP^r_d(\bar{t}, 1)$}
\obj(1,5)[G]{$\vdots$}
\obj(5,5)[H]{$\vdots$}
\obj(1,6)[I]{$\bbP^1[N,(r+1)(d-2)+1)]$}
\obj(5,6)[J]{$\bbP^r_d(\bar{t},d- 1)$}
\obj(1,7)[K]{$\bbP^1[N] = \bbP^1(N,N-1) $}
\obj(5,7)[L]{$G(\bar{t})$}
\obj(9,3)[M]{$\bbP^r_d$}
\obj(9,7)[N]{$\AG$}
\obj(11,5)[N1]{$\AM$}
\obj(9,6)[Z]{$\bbP^r_d(d-1)$}
\obj(9,4)[X]{$\bbP^r_d(1)$}
\obj(9,5)[W]{$\vdots$}
\obj(11,4)[Z1]{$\overline{M}_{0,0}(\bbP^r,d,d-1)$}
\obj(11,2)[X1]{$\overline{M}_{0,0}(\bbP^r,d,1)$}
\obj(11,3)[W1]{$\vdots$}
\mor{J}{Z}{}
\mor{F}{X}{}
\mor{N}{Z}{}
\mor{X}{M}{}
\mor{N1}{Z1}{}
\mor{X1}{B}{}
\mor{X}{X1}{}
\mor{Z}{Z1}{}
\mor{D}{C}{$(\bbC^*)^r$}[\atright,\solidarrow]
\mor{E}{C}{}
\mor{F}{E}{}
\mor{F}{D}{}
\mor{G}{E}{}
\mor{H}{F}{}
\mor{K}{I}{}
\mor{L}{J}{}
\mor{L}{K}{$(\bbC^*)^r$}[\atright,\solidarrow]
\mor{J}{I}{}
\mor{I}{G}{}
\mor{J}{H}{}
\mor{L}{N}{$(S_d)^{r+1}$}
\mor{D}{M}{$(S_d)^{r+1}$}
\mor{N}{N1}{}
\mor{M}{B}{}
\mor{W}{X}{}
\mor{Z}{W}{}
\mor{W1}{X1}{}
\mor{Z1}{W1}{}

\enddc
\end{figure}

Notice that $\bbP^r_d / G = \overline{M}_{0,0}(\bbP^r,d,1) = \dots = \overline{M}_{0,0}(\bbP^r,d,\frac{d+1}{2})$.  This is because up to that point, the exceptional loci of the blow ups will lie outside the stable locus.  For example, the exceptional divisor of $\bbP^r_d(\bar{t},1)\to \bbP^r_d(\bar{t})$ corresponds to a curve with two components.  One component is parametrized, and the other has all the $d(r+1)$ points on it.


\bibliographystyle{plain} 
\bibliography{ParkerStablemap}        

\begin{thebibliography}{10}

\bibitem{BF}
Gilberto Bini and Claudio Fontanari.
\newblock On the cohomology of {$\overline{M}_{0,n}(\bbP^1,d)$}.
\newblock {\em Commun. Contemp. Math.}, 4(4):751--761, 2002.

\bibitem{D}
Igor Dolgachev.
\newblock {\em Lectures on invariant theory}, volume 296 of {\em London
  Mathematical Society Lecture Note Series}.
\newblock Cambridge University Press, Cambridge, 2003.

\bibitem{FP}
W.~Fulton and R.~Pandharipande.
\newblock Notes on stable maps and quantum cohomology.
\newblock In {\em Algebraic geometry---Santa Cruz 1995}, volume~62 of {\em
  Proc. Sympos. Pure Math.}, pages 45--96. Amer. Math. Soc., Providence, RI,
  1997.

\bibitem{G}
Alexander~B. Givental.
\newblock Equivariant {G}romov-{W}itten invariants.
\newblock {\em Internat. Math. Res. Notices}, (13):613--663, 1996.

\bibitem{hart}
Robin Hartshorne.
\newblock {\em Algebraic geometry}.
\newblock Springer-Verlag, New York, 1977.
\newblock Graduate Texts in Mathematics, No. 52.

\bibitem{H1}
Yi~Hu.
\newblock Relative geometric invariant theory and universal moduli spaces.
\newblock {\em Internat. J. Math.}, 7(2):151--181, 1996.

\bibitem{H2}
Yi~Hu.
\newblock Moduli spaces of stable polygons and symplectic structures on
  {$\overline{M}{}\sb {0,n}$}.
\newblock {\em Compositio Math.}, 118(2):159--187, 1999.

\bibitem{HK}
Yi~Hu and Sean Keel.
\newblock Mori dream spaces and {GIT}.
\newblock {\em Michigan Math. J.}, 48:331--348, 2000.
\newblock Dedicated to William Fulton on the occasion of his 60th birthday.

\bibitem{K}
Sean Keel.
\newblock Intersection theory of moduli space of stable {$n$}-pointed curves of
  genus zero.
\newblock {\em Trans. Amer. Math. Soc.}, 330(2):545--574, 1992.

\bibitem{k1}
Frances~Clare Kirwan.
\newblock {\em Cohomology of quotients in symplectic and algebraic geometry},
  volume~31 of {\em Mathematical Notes}.
\newblock Princeton University Press, Princeton, NJ, 1984.

\bibitem{mori}
J{\'a}nos Koll{\'a}r and Shigefumi Mori.
\newblock {\em Birational geometry of algebraic varieties}, volume 134 of {\em
  Cambridge Tracts in Mathematics}.
\newblock Cambridge University Press, Cambridge, 1998.
\newblock With the collaboration of C. H. Clemens and A. Corti, Translated from
  the 1998 Japanese original.

\bibitem{MFK}
D.~Mumford, J.~Fogarty, and F.~Kirwan.
\newblock {\em Geometric invariant theory}, volume~34 of {\em Ergebnisse der
  Mathematik und ihrer Grenzgebiete (2) [Results in Mathematics and Related
  Areas (2)]}.
\newblock Springer-Verlag, Berlin, third edition, 1994.

\bibitem{M}
Andrei Musta\c{t}\v{a} and Anca Musta\c{t}\v{a}.
\newblock Intermediate moduli spaces of stable maps.
\newblock ARXIV: math.AG/0409569.

\bibitem{P1}
Rahul Pandharipande.
\newblock The {C}how ring of the nonlinear {G}rassmannian.
\newblock {\em J. Algebraic Geom.}, 7(1):123--140, 1998.

\end{thebibliography}
\index{Bibliography@\emph{Bibliography}}
\end{document}